\documentclass[onepage,11pt,a4paper,pdftex]{article}

\usepackage{titletoc}
\usepackage{graphicx, amsfonts, amsmath, amssymb, amscd,  theorem, exscale,
epic, eepic, epsfig}

\usepackage{hyperref} %

\usepackage{makeidx}
\makeindex

\newcounter{Figure}
\theoremstyle{plain}
\theoremheaderfont{\scshape}

\newtheorem{The}{\bf Theorem}

\theorembodyfont{\upshape}

\theoremheaderfont{\scshape\small} \theorembodyfont{\scshape\small}

\newcommand{\slap}{\mbox{$ \triangle \mkern -13mu / \ $}}
\newcommand{\nlap}{\mbox{$ \nabla \mkern -13mu / \ $}}
\newcommand{\dlap}{\mbox{$ div \mkern -13mu / \ $}}

\newcommand{\be}{\begin{equation}}
\newcommand{\ee}{\end{equation}}
\newcommand{\bea}{\begin{eqnarray}}
\newcommand{\eea}{\end{eqnarray}}
\newcommand{\beas}{\begin{eqnarray*}}
\newcommand{\eeas}{\end{eqnarray*}}

\newcommand{\sun}{\mbox{$ \circ \mkern -7mu   \cdot \  $}}

\begin{document}

\begin{center}
{\Large \bf  The Electromagnetic Christodoulou Memory Effect \\
\vspace{5pt}
in Neutron Star Binary Mergers
} \\ 
\end{center}
\vspace{5pt}
\begin{center}
\footnote{L. Bieri is supported by NSF grant DMS-0904583 
and S.-T. Yau is supported by NSF
grant PHY-0937443 and DMS-0904583. \\ 
Lydia Bieri, University of Michigan, Department of Mathematics, Ann Arbor MI.
lbieri@umich.edu \\ 
PoNing Chen, Harvard University, Department of Mathematics, Cambridge MA.
pchen@math.harvard.edu \\ 
Shing-Tung Yau, Harvard University, Department of Mathematics, Cambridge MA.
yau@math.harvard.edu}
%
{\large \bf Lydia Bieri, PoNing Chen, Shing-Tung Yau } \\ 
\end{center}
\vspace{5pt}
\noindent
\begin{abstract} 
Gravitational waves are predicted by the general theory of relativity. 
In \cite{chrmemory} D. Christodoulou showed that gravitational waves have a
nonlinear memory. 
We proved in \cite{lpst1} that the electromagnetic field contributes at highest
order to the nonlinear memory effect of gravitational waves. 
In the present paper, we study this electromagnetic Christodoulou memory effect and
compute it for binary neutron star mergers. These are typical sources of
gravitational radiation. 
During these processes, not only mass and momenta are radiated away in form of 
gravitational waves, but also very strong magnetic fields are produced and
radiated away. 
Thus the observed effect on test masses of a laser interferometer gravitational wave
detector will be enlarged by the contribution 
of the electromagnetic field. Therefore, the present results are important for the
planned experiments. 
Looking at the null asymptotics of spacetimes, which are solutions of the
Einstein-Maxwell (EM) equations, we derived in \cite{lpst1} the electromagnetic
Christodoulou memory effect. 
Moreover, our results allow to answer astrophysical questions, as 
the knowledge about the amount of energy radiated away in a neutron star binary
merger enables us to gain information about the source of the gravitational waves. 
\end{abstract}
The main goal of this paper is to discuss the electromagnetic Christodoulou memory
effect of gravitational waves and to compute 
this effect for typical sources. 
In \cite{chrmemory} D. Christodoulou showed that gravitational waves have a
nonlinear memory. 
In our paper \cite{lpst1} we proved 
that for spacetimes solving the Einstein-Maxwell (EM) equations, the
electromagnetic field 
contributes at highest order to the nonlinear memory effect of gravitational waves. 
In the present paper, we also calculate it for neutron star binary mergers. We find
that for typical 
constellations, very strong magnetic fields enlarge this effect considerably. 
Fields which are strong enough have so far only been known to be produced during
mergers of neutron star binaries. 
The latter are well known to be frequent events. There is a vast astrophysical
literature about this. \\ \\ 
Moreover, our results in \cite{lpst1} and in the present paper, are also important
from a purely astrophysical point of view. 
Namely, the knowledge about the amount of energy radiated away in a neutron star
binary merger allows to tell 
in the experiment what type of source the gravitational waves are coming from. 
Thus, our findings in the gravitational wave experiment will contribute to
astrophysical results. \\ \\ 
A major goal of general relativity (GR) and astrophysics is to precisely describe
and finally observe gravitational radiation, one of the 
predictions of GR. We know from the work \cite{chrmemory} of Christodoulou that also
these waves radiate. 
That is, in a laser interferometer gravitational wave detector, this will show in a
permanent displacement 
of test masses after a wave train passed. 
The latter is known as the 
Christodoulou nonlinear memory effect. 
In \cite{chrmemory} Christodoulou showed how the nonlinear memory effect 
can be measured as a permanent displacement of test
masses in such a detector. He derived
a precise formula for this permanent displacement in the Einstein vacuum (EV) case. 
The present authors proved in \cite{lpst1} that when electromagnetic fields are
present, they will contribute to this nonlinear effect at 
highest order. 
In fact, we showed that 
for the EM equations this permanent 
displacement exhibits a term coming from the electromagnetic field, which is at
the same highest order as the purely gravitational term that governs the EV
situation. Moreover, we showed that the instantaneous displacement
of the test masses is not changed at leading order by the electromagnetic field. 
To see this, we investigated spacetimes of solutions of the Einstein-Maxwell (EM)
equations at null infinity. \\ \\ 
Typical sources for gravitational waves are binary neutron star mergers and binary
black hole mergers. As the former are known to be much more frequent, it is likely
that gravitational waves as well as the nonlinear memory effect will first be
measured from binary neutron star mergers. During such processes mass and momenta
are radiated away. Moreover, large magnetic fields are produced and radiated
away. The radiation travels at the speed of light. 
That means, it moves along null hypersurfaces of corresponding spacetimes. 
Therefore, in order to fully understand all the different situations, one has to
investigate spacetimes which are solutions of the 
Einstein equations. Taking into account the strong magnetic fields which are
generated during binary neutron star mergers, we consider spacetimes solving the
Einstein-Maxwell equations. 
As the sources are very far away, we can think of us as doing the experiment at null
infinity. Therefore it is very important to understand the geometry of spacetimes
especially at null infinity, that is when we let $t \to \infty$ along null
hypersurfaces in the corresponding spacetimes.  \\ \\ 
In this paper, we discuss the electromagnetic Christodoulou memory effect and
compute concrete examples for binary neutron star mergers. 
In \cite{lpst1}, we derived this effect in the regime of the EM equations. 
First, we recall the Bondi mass
loss formula obtained in \cite{zip2} for spacetimes solving the EM equations. 
\begin{equation} \label{bondimassloss}
\frac{\partial }{\partial u}M\left( u\right) =\frac{1}{8\pi }%
\int_{S^{2}}\left( \left| \Xi \right| ^{2}+\frac{1}{2}\left| A_{F}\right|
^{2}\right) d\mu _{\overset{\circ }{\gamma }}
\end{equation}
Compared to the formula obtained in \cite{sta} for spacetimes solving the EV
equations, we have an additional term, $|A_F|^2$, from the electromagnetic field.
(See \cite{lpst1}.) \\  \\
As shown in the work of Christodoulou \cite{chrmemory}, $\Sigma^+ - \Sigma^-$ is the
term which governs the permanent displacement of test particles. 
Using this fact, 
Christodoulou shows that the gravitational field has a non-linear ``memory" which
can be detected by a gravitational-wave experiment in a spacetime solving the EV
equations. 
Here, $\Sigma$ denotes the asymptotic shear of outgoing null hypersurfaces $C_u$
that are level sets of a foliation by an 
optical function $u$, which we will discuss below. $\Sigma^+$ and $\Sigma^-$ are the
limits of $\Sigma$ as 
$u$ tends to $+ \infty$ respectively $- \infty$. \\ \\ 
In our paper \cite{lpst1}, we study the permanent displacement formula
for uncharged test particles of the same gravitational-wave experiment in a
spacetime solving the EM equations. 
We derive $\Sigma^+ - \Sigma^-$ in the
EM case, and we find that the electromagnetic field
changes the leading order term of the permanent displacement of test particles. 
Moreover, investigating the experiment for our setting in \cite{lpst1}, we prove
that the electromagnetic field does not enter the leading
order term of the Jacobi equation. As a result, to leading order, it does not change
the instantaneous displacement of test particles. 
But the electromagnetic field does contribute at highest order to the nonlinear
effect of the permanent displacement of test masses. \\ \\ 
To study the effect of gravitational waves, 
we follow the method introduced  by Christodoulou in \cite{chrmemory}. 
The analysis is based on the asymptotic behavior of
the gravitational field obtained at null and spatial infinity. 
These rigorous
asymptotics allow us to study the structure of the spacetimes at null infinity. 
To foliate the spacetime, we use a 
time function $t$ and an optical function $u$. 
We denote the corresponding lapse functions by $\phi$ respectively $a$.  
Whereas each level set of
$t$, $H_t$ is a maximal spacelike hypersurface, 
each level set of $u$, $C_u$, is an outgoing null hypersurface.  
Along the null hypersurface $C_u$, we pick a suitable pair of normal
vectors. 
The flow along these vector fields generates a
family of diffeomorphisms $\phi_u$ of $S^2$. 
Using $\phi_u$ we pull back tensor
fields in our spacetime. 
In this manner, we can study their limit at null infinity along
the null hypersurface $C_u$.  
Building on these, we then take the limit as $u$ goes to $\pm \infty$, which allows
us to 
investigate the effect of gravitational waves. 
For a detailed explanation of the structure at null infinity, see \cite{chrmemory}
by Christodoulou. \\ \\
Understanding gravitational radiation and therefore null infinity heavily relies on
the rigorous understanding of 
the corresponding spacetimes. 
The methods introduced in \cite{sta}, used in \cite{zip}, \cite{zip2} and
\cite{lydia1}, \cite{lydia2}, reveal the
structure of 
the null asymptotics of our spacetimes. 
In these works, stability results were proven. The authors showed that under a
smallness condition on asymptotically flat initial data for the EV respectively EM
equations, this can be extended uniquely to a smooth, globally hyperbolic and
geodesically complete spacetime solving the EV respectively EM equations. The
spacetime obtained is globally asymptotically flat. 
The main achievements are generally two-fold: First, existence and uniqueness
theorems were proven. To ensure these, one has to impose smallness conditions.
Second, precise descriptions of the asymptotic behavior of the spacetimes were
derived. We stress the fact, that the results about null infinity are largely
independent of the smallness. 
An elaborated geometric-analytic procedure led to these results. And many
mathematical theorems were proven on the way. 
However, the outcome exhibits a physical result in point two from which
Christodoulou in \cite{chrmemory} derived the Christodoulou memory effect of
gravitational waves in
the EV case and the present authors 
in \cite{lpst1} in the EM case. 
In what follows, let us, discuss the new physical results and compute the effects
for different binary neutron star constellations. \\ \\ 
First we recall the Einstein-Maxwell equation. The electromagnetic field is
represented by a 
skew-symmetric 2-tensor $F_{\mu\nu}$. The stress-energy tensor corresponding to
$F_{\mu\nu}$ is 
\[ T_{\mu \nu} = \frac{1}{4 \pi} \big(  F_{\mu}^{\ \rho} F_{\nu \rho} -
\frac{1}{4}g_{\mu \nu} F_{\rho \sigma} F^{\rho \sigma}  \big) \]
The Einstein-Maxwell equations read:
\begin{equation}
 \label{EM}
\begin{split}
\ R_{\mu \nu}  \ =& \ 8 \pi T_{\mu \nu}  \\
D^{\alpha} F_{\alpha \beta}  = & 0 \\
D^{\alpha} \ ^* F_{\alpha \beta} = & 0. 
\end{split}
\end{equation}
Let $S_{t,u}$ be the intersection of the hypersurface $H_t$ and the null cone $C_u$.
 Let $N$ be the spacelike unit normal vector of $S_{t,u}$ in $H_t$ and $T$ be the 
timelike unit normal vector of $H_t$ in the spacetime. Let $\{ e_a \}_{a=1,2}$  be
an orthonormal frame on $S_{t,u}$. We have the following orthogonal frame $(T, N,
e_2, e_1)$. 
This also gives us a pair of null normal vectors to $S_{t,u}$, namely $L=T+N$ and
$\underline{L}=T-N$. Together with $\{ e_a \}_{a=1,2}$, they 
form a null frame. The following is a picture of the null cone $C_u$ together with
the null frame $(L, \underline{L}, e_2, e_1)$.  \\ \\ 
\begin{center}
\includegraphics[scale=0.45]{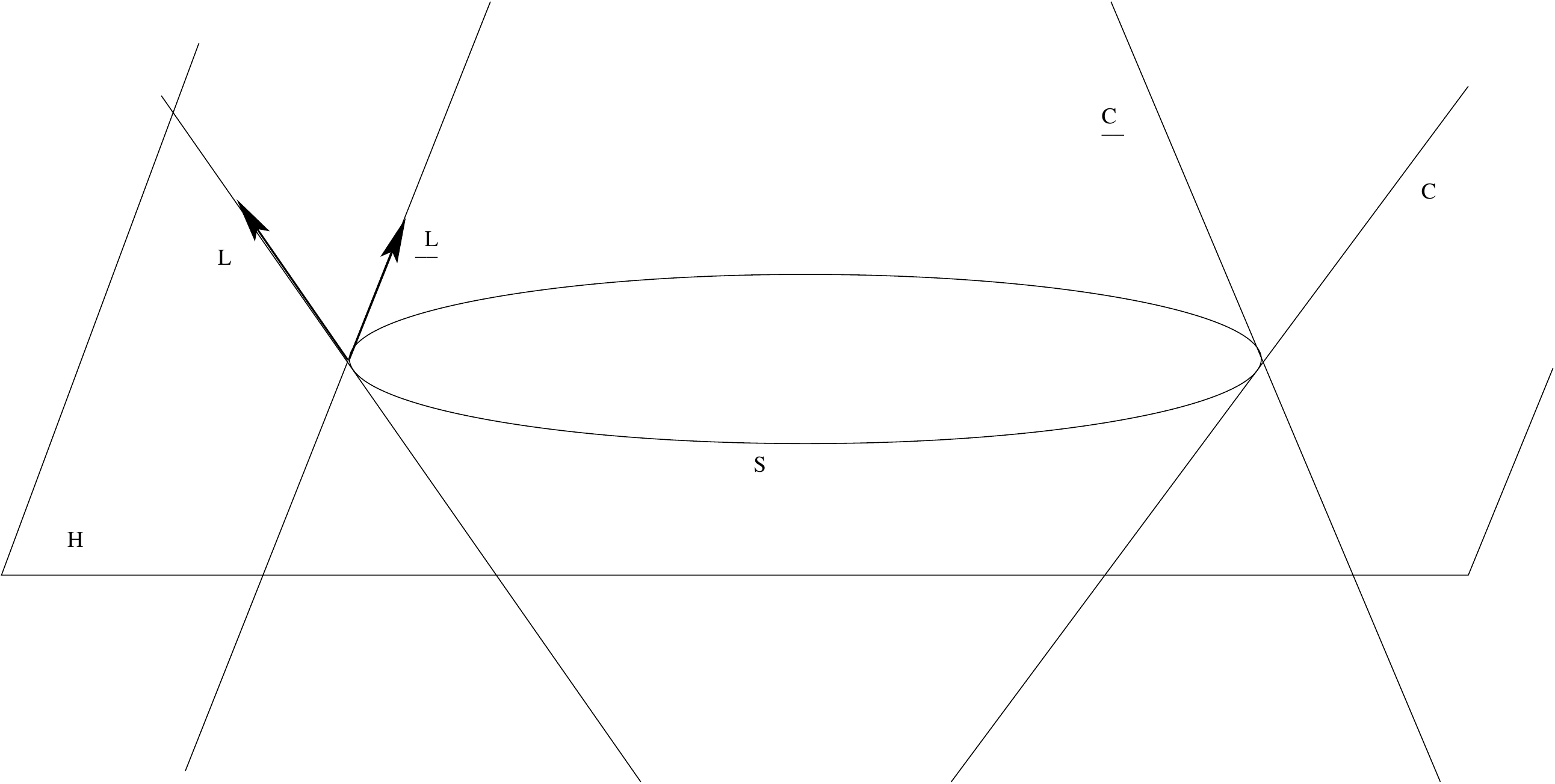}
\end{center}
We can decompose the Weyl curvature tensor and the electromagnetic field with
respect to the null frame or the orthogonal frame. The asymptotics of these
components are studied in 
\cite{zip} and \cite{zip2}. These asymptotics are important for the understanding of
the geometry of null infinity. For simplicity, we will only list
the components of the spacetime curvature and  electromagnetic field that are used
in our discussion. Please see 
\cite{zip} and \cite{zip2} for more details on the asymptotics. \\ \\ 
Let $X, Y$ be arbitrary tangent vectors to $S$ at a point in $S$. 
Given the null frame $e_4 = L $, $e_3 = \underline{L}$ and  $\{ e_a \}_{a=1,2}$, let 
$\chi (X, Y) = g(\nabla_X L, Y)$ 
and $\underline{\chi} (X, Y) = g(\nabla_X \underline{L}, Y)$ be the second
fundamental forms with respect to $L$
and $\underline{L}$, respectively. Let $\widehat \chi$ and $ \widehat{
\underline{\chi} }$ be their traceless parts.  
We also need the following null components of the Weyl curvature
\[ \underline \alpha_W(X,Y) = R(X, \underline L, Y , \underline L) \]
and the electromagnetic field 
\begin{equation}
\begin{array}{lll}
F_{A3}=\underline{\alpha }(F)_A  &  & F_{A4}=\alpha \left( F\right)_A
\\ 
F_{34}=2\rho \left( F\right) &  & F_{12}=\sigma \left( F\right) 
\end{array}
\label{null-electricfield}
\end{equation}
We have the following limit of the above quantities at null infinity
\[ \lim_{C_u , t \to \infty} r^2 \widehat \chi = \Sigma   \text{\, , \,}  \lim_{C_u
, t \to \infty} r \widehat{ \underline{\chi} }= 2 \Xi   \]  
\[ \lim_{C_u , t \to \infty} r \underline \alpha_W = A_W    \text{\, , \,}  \lim_{C_u , t \to \infty} r \underline{\alpha}_F = A_F     \]
As shown in \cite{chrmemory}, the permanent displacement of the test masses of a
laser interferometer gravitational-wave detector is 
governed by $\Sigma^+ - \Sigma^-$ where 
\begin{equation*} \lim_{u \to \pm \infty} \Sigma =  \Sigma ^{\pm} \end{equation*}
\begin{The} \label{theXiSigmaA}  \cite{zip}, \cite{zip2}
We have the following equations for $\Sigma$, $\Xi$ and $A_W$
\[ \frac{\partial \Sigma}{\partial u} = - \Xi  \text{ \, and \,} \frac{\partial \Xi}{\partial u} = - \frac{1}{4} A_W \]
\end{The}
In our paper \cite{lpst1}, we prove that 
in a spacetime solving the Einstein--Maxwell equations, \\ 
$\Sigma^+ - \Sigma^-$ is governed by the following relation. 
\begin{The} \label{displ*1} \cite{lpst1}
Let 
\be \label{Thm*FXiAF*1}
F (\cdot)  =   \int_{- \infty}^{\infty} 
\big( 
\mid \Xi (u, \cdot) \mid^2 + \frac{1}{2} \mid A_F (u, \cdot) \mid^2  
\big)
du  \ \ . 
\ee
Then 
$\Sigma^+ - \Sigma^-$ is given by the following equation on $S^2$: 
\be \label{Thm*divSigma+-*2}
\stackrel{\circ}{\dlap} (\Sigma^+ - \Sigma^-) = \stackrel{\circ}{\nlap} \Phi \ \ . 
\ee
where $\Phi$ is the solution with $\bar{\Phi} = 0$ on $S^2$ of the equation 
\[
\stackrel{\circ}{\slap} \Phi = F - \bar{F}   \ \ . 
\]
\end{The}
Comparing this with the EV case studied in the last chapter of \cite{sta} and used
in \cite{chrmemory}, 
where the corresponding formula was 
$F (\cdot)  =   \int_{- \infty}^{\infty} 
\mid \Xi (u, \cdot) \mid^2 du$, we find that new the electromagnetic part
$\frac{1}{2} \mid A_F (u, \cdot) \mid^2$ appears in the integral. 
In fact, in our proof, we derive the limiting formulas and obtain the said
electromagnetic contribution in $\Sigma^+ - \Sigma^-$. (See \cite{lpst1}.)  \\ \\
{\bf Gravitational Wave Experiment} \\ \\ 
How will our findings relate to experiment? In what follows, we are going to show
how the 
electromagnetic field enters the experiment. In particular, we will discuss the
instantaneous and the permanent displacement of test masses. 
For a detailed explanation of the experiment we refer to \cite{chrmemory} and for a
detailed derivation in the EM case we refer to \cite{lpst1}. \\ \\ 
Consider a laser interferometer gravitational-wave
detector with three test masses. 
We denote the reference mass by $m_0$, this is also the location of the beam
splitter. The masses $m_0$, $m_1$, $m_2$ are 
suspended by equal length pendulums of length $d_0$. The motion of the masses in the
horizontal plane can be considered free for timelike scales much 
shorter than the period of the pendulums. 
Now one measures the distance of the masses $m_1$ and $m_2$ from the reference test
mass $m_0$ by laser interferometry. 
We observe a difference of phase of the laser light at $m_0$ whenever the light
travel times between $m_0$ and $m_1$, $m_2$, respectively, differ. \\ \\ 
The motion of the masses $m_0$, $m_1$, $m_2$ is described by geodesics $\gamma_0$,
$\gamma_1$, $\gamma_2$ in spacetime. 
Denote by $T$ the unit future-directed tangent vectorfield of $\gamma_0$ and by $t$
the arch length along $\gamma_0$. 
Let then $H_t$ be for each $t$ the spacelike, geodesic hyperplane through $\gamma_0
(t)$ orthogonal to $T$. 
At $\gamma_0 (0)$ pick an orthonormal frame $(E_1, E_2, E_3)$ for $H_0$. By
parallelly propagating it along $\gamma_0$, 
we obtain the orthonormal frame field 
$(T, E_1, E_2, E_3)$ along $\gamma_0$, where at each $t$ the $(E_1, E_2, E_3)$ is an
orthonormal frame for $H_t$ at $\gamma_0 (t)$. 
Then we can assign to a point $p$ in spacetime close to $\gamma_0$ and lying in
$H_t$ the cylindrical normal coordinates $(t, x^1, x^2, x^3)$.  \\ \\  
Supoose that the distance $d$ is much smaller than the time scale in which the
curvature of the spacetime varies significantly. Then 
the geodesic deviation from $\gamma_0$, namely the Jacobi equation (\ref{jacobi*1}),
replaces the geodesic equation for $\gamma_1$ and $\gamma_2$. 
Let $R_{k0l0} = R(E_k, T, E_l, T)$, then we write 
\be \label{jacobi*1}
 \frac{d^2 x^k}{dt^2}=-R_{k0l0} x^l 
 \ee
We can decompose $R_{k0l0}$ into the Weyl curvature and the Ricci curvature
\[ R_{k0l0}=W_{k0l0} +
\frac{1}{2}(g_{kl}R_{00}+g_{00}R_{kl}-g_{0l}R_{k0}-g_{0k}R_{l0}). \]
From the EM equations (\ref{EM}) we find 
\be \label{R00}
R_{00} =  \frac{1}{2} ( \mid \underline{\alpha} (F) \mid^2 + \mid \alpha (F) \mid^2 
)  + \rho (F)^2 + \sigma (F)^2 
\ee 
The component $R_{00}$ observes the term $\mid \underline{\alpha} (F) \mid^2$, where
$\underline{\alpha} (F)$ is the electromagnetic field component with worst decay
behavior, but entering $R_{00}$ as a quadratic. 
Hence, $R_{00}$ is of the order $O(r^{-2})$. Whereas the leading order component of
the Weyl curvature is of the order $O(r^{-1})$. 
We give a detailed proof in our paper \cite{lpst1}. 
Thus, the electromagnetic field does not contribute at highest order to the
deviation measured by the Jacobi equation. 
As a consequence, it does only change at lower order the instantaneous displacement
of the test masses. 
However, we are going to see that it does change the nonlinear memory effect. \\ \\ 
Using the relations from theorem \ref{theXiSigmaA} and our theorem \ref{displ*1} as
well as the fact that $\Xi \to 0$ for $u \to \infty$ and taking the limit 
$t \to \infty$, we conclude that the test masses experience permanent displacements
after the passage of a wave train. In particular, 
this overall displacement of the test masses is described by $\Sigma^+ - \Sigma^-$ 
\be \label{overall} 
\Delta x^A_{(B)} = - \frac{d_0}{r}(\Sigma^+_{AB}-\Sigma^-_{AB}) 
\ee
where from our theorem \ref{displ*1} one sees that the right hand side of
(\ref{overall}) includes the electromagnetic field terms at highest order. \\ \\ 
Let us now derive formula (\ref{overall}). We will use $L = T - E_3$ and
$\underline{L} = T + E_3$. 
Then we write the leading components of the curvature $\underline{\alpha}_{AB} (W)$
and of the 
electromagnetic field $\underline{\alpha}_A (F)$ as follows: 
\beas
\underline{\alpha}_{AB} (W) \ & = & \ R \ (E_A, \ \underline{L}, \ E_B, \
\underline{L}) \ = \  \frac{A_{AB}(W)}{r} \ + \ o \ (r^{-2})  \\ 
\underline{\alpha}_{A} (F) \ & = & \ F (E_A, \underline{L})   \ = \ 
\frac{A_{A}(F)}{r} \ + \ o \ (r^{-2}) 
\eeas
Let $x^k_{(A)}$ with $A = 1,2$ denote the $k^{th}$ Cartesian coordinate of the mass
$m_A$. 
From \cite{chrmemory} and \cite{lpst1} one sees that there is no acceleration to
leading order in the vertical direction. 
One starts with $m_1$, $m_2$ being at rest at equal distance $d_0$ from $m_0$ at
right angles from $m_0$. 
Thus to leading order it is 
\be 
\ddot{x}^A_{(B)} \ = \ - \frac{1}{4} r^{-1} d_0 A_{AB} 
\ee
In particular, the initial conditions are as $t \to - \infty$: \\ 
$x^B_{\ (A)} =  d_0  \delta^B_A  \ , \ 
\dot{x}^B_{\ (A)}  =  0 \ , \ 
x^3_{\ (A)}  =  0 \ , \ \dot{x}^3_{\ (A)}  = 0$. \\ \\ 
Integrating gives 
\be 
\dot{x}^A_{\ (B)} \ (t) \ = \ - \ \frac{1}{4} \  d_0 \  r^{-1} \ 
\int_{- \infty}^t \ A_{AB} \ (u) \ d u \ . 
\ee
From theorem \ref{theXiSigmaA} equation $ \frac{\partial \Xi}{\partial u} = -
\frac{1}{4} A_W $ and 
$lim_{\mid u \mid \to \infty} \Xi = 0$, one substitutes and concludes 
\be
\dot{x}^A_{\ (B)} \ (t) \ = \ \frac{d_0}{r}  \ \Xi_{AB} \ (t) \ . 
\ee
As $\Xi \to 0$ for $u \to \infty$, the test masses return to rest after 
the passage of the gravitational waves. 
Now, we use theorem \ref{theXiSigmaA} equation $\frac{\partial \Sigma}{\partial u} =
- \Xi$ and 
integrate again to obtain 
\be
x^A_{\ (B)} \ (t) \ = \ - \ (\frac{d_0}{r}) \ (\Sigma_{AB} \ (t) \ - \ \Sigma^-) \ . 
\ee
Finally, by taking the limit $t \to \infty$ one derives that the test masses obey
permanent displacements. 
This means that 
$\Sigma^+ - \Sigma^-$ is equivalent to an overall displacement of the test masses
given by (\ref{overall}):
\[
\triangle \ x^A_{\ (B)} \ = \ - \ (\frac{d_0}{r}) \ (\Sigma^+_{AB} \ - \
\Sigma^-_{AB}) \ . 
\]
The right hand side of (\ref{overall}) includes terms from the
electromagnetic field at highest order as given in our theorem \ref{displ*1}. \\ \\ 
In the next subsection, we are going to apply our results to astrophysical data for
binary neutron star mergers. \\ \\ 
{\bf Binary neutron star mergers} \\ \\ 
We compute the electromagnetic Christodoulou memory effect for typical sources, that
is for different constellations of binary neutron star (BNS) mergers. \\ \\ 
In a binary neutron star or binary black hole system, the two objects are orbiting
each other. In Newtonian physics, they would stay like that forever. However,
according to the theory of general relativity such a system must radiate away
energy. Therefore, the radius of the orbits must shrink and finally the objects will
merge. \\ \\ 
As binary neutron star systems are much more frequent than binary black hole
systems, it is very likely that gravitational waves as well as the nonlinear memory
effect of gravitational waves will be 
detected first from the former systems. 
The magnetic fields produced and radiated away during the merger of two neutron
stars are among the largest magnetic fields 
known in astrophysics. In fact, in the electromagnetic Christodoulou memory effect
that we derived, the magnetic field enlarges the 
nonlinear displacement of (non-charged) test masses significantly. As we are going
to show in this subsection, the contribution from the 
magnetic field is very important, as it is very big for a large part of the known
constellations. \\ \\ 
Astrophysical data gives for typical neutron star binaries a range of possible
constellations which allow the mass and the 
magnetic field to vary within given boundaries. Typically, the mass of a neutron
star is around slightly more than $1 M_{\sun}$ and the radius 
of a neutron star is $3$ - $ 30$ km. Thus, the typical mass for a BNS system ranges
between $2.6$ and $2.8 M_{\sun}$ \footnote[1]{$ M_{\sun} = 1$ solar mass $ \approx
1.9891 \cdot 10^{33}$ g}. 
In such a system, as the neutron stars are spiraling 
around each other, they are radiating away gravitational and magnetic energy. 
The inspiral goes with increasing speed and the BNS system emits an increasing
amount of electromagnetic and gravitational energy, which becomes extremely large
when the orbit radius is about $10$ - $100 $km.  
For the detection of the electromagnetic Christodoulou effect, the largest
contribution will come from the last phase of the inspiral, starting when the orbit
radius is about 10 times the neutron star radius. 
In the literature, we find that the merger times range from a few milli-seconds up
to $1000$ ms. We would like to compare the amount of 
gravitational energy radiated away during the merger to the amount of magnetic
energy radiated away. On the one hand, the amount 
of gravitational energy radiated away is well-known. In general, about $1 \%$ of the
initial mass is radiated away during a merger. This is about 
$10^{52}$ erg.\footnote[2]{$1 \  {\rm erg} = 1 {\rm g \cdot cm^2s^{-2}}$ and $1 M_{\sun}
\approx 1.78 \cdot 10^{54}$ erg} On the other hand, the amount of magnetic energy radiated
away could vary 
drastically depending on different constellations. Typically, the rate of change for
the magnetic field is $\frac{dB}{dt} \approx 10^9 - 10^{17} \,{\rm G(ms)^{-1}}$ and the
magnetic field produced in 
the merger is about $10^9 - 10^{17} \,{\rm G}. \footnote[3]{$1$ Gauss: $1 G = 10^{-4}\, {\rm kg \cdot
C^{-1} s^{-1}} = 10^{-1}\,{\rm g \cdot C^{-1}  s^{-1}}$  }$ \\ \\ 
Comparing the energy from the radiated mass, i.e. purely gravitational, and from the 
magnetic field, we observe that during the merger of BNS very large
magnetic fields are produced and radiated away in certain scenarios. \\ \\ 
Consider the following data.
Assume: Total mass of BNS is initially $ 2 M_{\sun}$, $1 \%$ of the total mass will
be radiated away during the (whole) merger, radius of each neutron star is $10$ km.
Under the assumption, the gravitational energy radiated away is about $3.56 \times 10^{52} $\,erg.
\\ \\ 
In the physics literature, one finds many linearized models. However, the Einstein
equations being nonlinear, the main information usually gets lost in linearized
models. As we do investigate the nonlinear problem here, and as 
the results of \cite{chrmemory} and \cite{lpst1} show the Christodoulou memory
effect of gravitational waves to be a nonlinear phenomenon, we consider a
corresponding nonlinear model for the neutron star binary mergers. 
Thus, we use the results of Zipser's global stability work \cite{zip} and
\cite{zip2} for the 
initial value problem in 
spacetimes satisfying the Einstein-Maxwell equations. 
We assume that outside the neutron star, the magnetic field decays like $r^{-5/2}$.
Such decay at spatial infinity is suggested by the decay obtained in \cite{zip},
\cite{zip2}. 
One might want to consider situations with a slightly different decay of the
magnetic field. This would not affect the main picture, as one finds during the
computations that the decay of the magnetic field does not play a role here. 
Thus, we work with the nonlinear model explained in the following paragraph. \\ \\ 
Now, consider such a BNS system with 
the magnetic field $B$ 
initially being $B = 10^{13}$ G and $\frac{dB}{dt} = 10^{13}\, {\rm G(ms)^{-1}}$ on the
surface of the neutron star. 
Assume that the merger time is $1000$ ms. We estimate the total magnetic energy
radiated away using the following model. We assume that through the merger, the
matter of the neutron star
stays in a ball of radius $10$ km. We compute the contribution from the 
magnetic field outside 
the support of the matter of the neutron star. As a result, we simply use the vacuum
magnetic constant when computing
the magnetic energy density. 
Moreover, we assume that outside the neutron star, the magnetic field decays at the rate of 
$r^{-5/2}$.
Using this model, the energy radiated away from the 
magnetic field is about $4.78 \cdot 10^{49}$ erg. In this case, the addition
of a magnetic field has a small contribution to the memory effect. 
\\ \\ 
Next, consider a BNS system with the above data, but where the magnetic field $B$ is initially $B = 10^{15}$\,G and $\frac{dB}{dt} = 10^{15}\, {\rm G(ms)^{-1}}$ on the surface of the neutron star. Assume that the merger time is $1000$ ms. We compute that the total magnetic energy radiated 
away is about $4.78 \cdot 10^{53}$ erg. This will be one order of magnitude higher than the
gravitational energy radiated away. This situation is consistent with astrophysical data.  
Also, in the numeric simulation in \cite{grb1}, \cite{grb2} and \cite{grb3}, it is observed that the magnetic field could increase by two orders of magnitude
during merger when one starts with magnetic fields around $10^9$ to $10^{12}$ G. When we start with a stronger magnetic field, the merger
would take longer and allow more time for the magnetic field to build up. Magnetic fields of similar initial strength are used in the simulation of \cite{lset}. 
Their simulations suggest that the addition of the magnetic field cause observable differences in the dynamics and gravitational waveforms. Moreover, it is noticed that the most important role of magnetic fields are on the long term evolution. This is similar to our conclusion from theorem \ref{theXiSigmaA} and \ref{displ*1}. Namely, the addition of a magnetic field does not change the system instantaneously but it does
contribute to the nonlinear long-term permanent change.\\ \\ 
To compare, note that the amount of energy emitted in a binary black hole merger 
is expected to be as follows: a binary black hole system with equal mass and no spin
would lose about $4\%$ of the mass during merger. 
However, BNS mergers occur more often than black hole mergers. \\ \\
{\bf Conclusions:} We find that among the variety of different constellations of BNS
systems there is a large part for which the magnetic field contributes to the
Christodoulou effect at the same highest order as the purely gravitational term.  \\
\\ 
{\bf Acknowledgment:} We thank Demetrios Christodoulou for fruitful discussions and
his interest in this work. \\ \\ \\

 
%
\vspace{10pt}
{\scshape Lydia Bieri \\ 
Department of Mathematics \\ 
University of Michigan \\ 
Ann Arbor, MI 48109, USA} \\ 
lbieri@umich.edu \\ \\ 
{\scshape PoNing Chen \\  
Department of Mathematics \\ 
Harvard University \\
Cambridge, MA 02138, USA} \\ 
pchen@math.harvard.edu \\ \\ 
{\scshape Shing-Tung Yau \\ 
Department of Mathematics \\ 
Harvard University \\
Cambridge, MA 02138, USA} \\ 
yau@math.harvard.edu

\begin{thebibliography}{99} 
\bibitem{lydia1} L. Bieri.  
        \begin{itshape} An Extension of the Stability Theorem of the Minkowski Space
in General Relativity. \end{itshape}
        ETH Zurich, Ph.D. thesis.  \textbf{17178}. 
        Zurich. (2007).  
\bibitem{lydia2} L. Bieri.  
        \begin{itshape} Extensions of the Stability Theorem of the Minkowski Space
in General Relativity. Solutions of the Einstein Vacuum Equations.
\end{itshape}
        AMS-IP. Studies in Advanced Mathematics. Cambridge. MA. (2009).   
\bibitem{lpst1} L. Bieri, P. Chen, S.-T. Yau. 
  \begin{itshape} Null Asymptotics of Solutions of the Einstein-Maxwell Equations in
General Relativity and 
  Gravitational Radiation.    \end{itshape} 
Submitted. (2010).  \\ 
http://arxiv.org/abs/1011.2267         
\bibitem{BBM} H. Bondi, M. G. J. van der Burg and A. W. K. Metzner.  
\begin{itshape} Gravitational Waves in General Relativity. VII. Waves from
Axi-Symmetric Isolated Systems. \end{itshape} 
Proc. Roy. Soc. A. \textbf{269} (1962).  21-52
\bibitem{vanderBurg}  M. G. J. van der Burg. 
        \begin{itshape} Gravitational Waves in General Relativity X. Asymptotic
Expansions for the Einstein-Maxwell Field \end{itshape}
         Proc. Roy. Soc. A. \textbf{310} (1969).  221-230
\bibitem{chrmemory}  D. Christodoulou. 
        \begin{itshape} Nonlinear Nature of Gravitation and Gravitational-Wave
Experiments. \end{itshape}
        Phys.Rev.Letters. \textbf{67}. 
        (1991). no.12. 1486-1489. 
\bibitem{chrdlbmathpgrt}  D. Christodoulou.        
        \begin{itshape} Mathematical problems of general relativity theory I and II.
\end{itshape}
         Volume 1: EMS publishing house ETH Z\"urich. (2008). Volume 2 to apppear:
EMS publishing house ETH Z\"urich. 
\bibitem{sta} D. Christodoulou, S. Klainerman.
        \begin{itshape} The global nonlinear stability of the Minkowski space.
\end{itshape}
        Princeton Math.Series \textbf{41}. 
        Princeton University Press. Princeton. NJ. (1993). 
\bibitem{grb1} Bruno Giacomazzo, Luciano Rezzolla, Luca Baiotti
 \begin{itshape} Can magnetic fields be detected during the inspiral of binary neutron stars?,
\end{itshape}
arXiv:0901.2722 [gr-qc]
\bibitem{grb2} Bruno Giacomazzo, Luciano Rezzolla, Luca Baiotti
 \begin{itshape} Accurate evolutions of inspiralling and magnetized neutron-stars:
equal-mass binaries,
\end{itshape}
arXiv:1009.2468v2 [gr-qc]
\bibitem{grb3} Luciano Rezzolla, Bruno Giacomazzo, Luca Baiotti, Jonathan Granot, Chryssa Kouveliotou, Miguel A. Aloy
 \begin{itshape} The missing link: Merging neutron stars naturally produce jet-like structures and can power short Gamma-Ray Bursts,
\end{itshape}
arXiv:1101.4298 [gr-qc]
\bibitem{lset} Yuk Tung Liu, Stuart L. Shapiro, Zachariah B. Etienne, Keisuke Taniguchi
 \begin{itshape} General relativistic simulations of magnetized binary neutron star
mergers,
\end{itshape}
arXiv:0803.4193v2 [astro-ph]
\bibitem{st}Shapiro, S. L. and Teukolsky, S. A.
\begin{itshape}
Black Holes, White Dwarfs and Neutron Stars: The Physics of Compact Objects
\end{itshape} 
Wiley-Vch New York, NY (1983)
\bibitem{zip} N. Zipser. 
        \begin{itshape} The Global Nonlinear Stability of the Trivial Solution of
the Einstein-Maxwell Equations.  \end{itshape}
        Ph.D. thesis. Harvard Univ. Cambridge MA. (2000).          
\bibitem{zip2} N. Zipser.  
        \begin{itshape} Extensions of the Stability Theorem of the Minkowski Space
in General Relativity. - Solutions of the Einstein-Maxwell Equations.
\end{itshape}
        AMS-IP. Studies in Advanced Mathematics. Cambridge. MA. (2009).           
\end{thebibliography}
\end{document}